  \documentclass[twoside,reqno,11pt]{fcaa}  %%% or in case of problems, use as below: %
\usepackage{graphicx}
\usepackage{epsfig}
\usepackage{amsthm}
\usepackage{amsmath}
\usepackage{latexsym}
\usepackage{amsfonts}
\usepackage{amssymb}

\newtheoremstyle{theorem}% name
  {15pt}          % space above
  {15pt}  % space below
  {\sl}  % bofy font
  {\parindent}
  {\sc}  % thm head font
  {. }   % punctuation after thm head
  { }    % space after thm head: `` ``=normal \newline=linebreak
  {}     % thm head specification
\theoremstyle{theorem}

\newtheoremstyle{defi}% name
  {15pt}          % space above
  {15pt}  % space below
  {\rm}  % bofy font
  {\parindent}     % ident - empty=no indent,  \parindent= paragraph indent
  {\sc}  % thm head font
  {. }    % punctuation after thm head
  { }    % space after thm head: `` ``=normal \newline=linebreak
  {}     % thm head specification
\theoremstyle{defi}

 \def\proofname{\indent\rm P\ r\ o\ o\ f.}
 \def\proofend{\hfill$\Box$}
 \def\qed{\hfill$\Box$} % instead of \square

 \usepackage{hyperref} % Editor will use to create hyperlinks %
 \def\theequation{\arabic{section}.\arabic{equation}}

 \def\themycopyrightfootnote{\vspace*{3pt}
 \copyright \, Year\,  Diogenes  Co., Sofia
   \par  \noindent pp. xxx--xxx, DOI: ......................
   \hfill  \vspace*{-36pt}
  }
 \title[ON DIFFERENTIABILITY OF SOLUTIONS \dots]
 {ON DIFFERENTIABILITY OF SOLUTIONS \\ OF FRACTIONAL DIFFERENTIAL EQUATIONS \\ WITH RESPECT TO INITIAL DATA \\ [3pt] IN ``FCAA'' JOURNAL}
 \author[\normalsize M. I. Gomoyunov]{\normalsize Mikhail I. Gomoyunov}

 \newtheorem{theo}{Theorem}[section]
 \newtheorem{lem}{Lemma}[section]
 \newtheorem{prop}{Proposition}[section]
 \newtheorem{cor}{Corollary}[section]

 \def\AC{\operatorname{AC}}
 \def\C{\operatorname{C}}
 \def\rd{\mathrm{d}}

\begin{document}

 \vbox to 2.5cm { \vfill }

%%% to make empty space of approx. 2.5cm %%%%%%
%%% will be replaced by Editor with the journal's and publoishers logos %%%%%%%%

 \bigskip \medskip

%%%% Abstract %%%%%%%%%%%%%%%%%%%%%%%%%
 \begin{abstract}

    In this paper, we deal with a Cauchy problem for a nonlinear fractional differential equation with the Caputo derivative of order $\alpha \in (0, 1)$.
    As initial data, we consider a pair consisting of an initial point, which does not necessarily coincide with the inferior limit of the fractional derivative, and a function that determines the values of a solution on the interval from this inferior limit to the initial point.
    We study differentiability properties of the functional associating initial data with the endpoint of the corresponding solution of the Cauchy problem.
    Stimulated by recent results on the dynamic programming principle and Hamilton--Jacobi--Bellman equations for fractional optimal control problems, we examine so-called fractional coinvariant derivatives of order $\alpha$ of this functional.
    We prove that these derivatives exist and give formulas for their calculation.

 \medskip


 \smallskip


 \end{abstract}

 \maketitle

%%%%%%% end make title %%%%%%%%%%%%%%%%%%%%%%%%%%%%%%%%%%
 \vspace*{-16pt}

%%%%%%%% begin papers' body %%%%%%%%%%%%%%%%%%%%%%%%%%%%%

%%%%%%%%%%%%%%%%%%%%%%%%%%% Section 1 %%%%%%%%%%%%%
\section{Introduction}
\setcounter{section}{1}
\setcounter{equation}{0}\setcounter{theorem}{0}

    In this paper, we deal with the Cauchy problem for the nonlinear fractional differential equation
    \begin{equation} \label{differential_equation_t}
        (^C D_{0 +}^\alpha x)(\tau)
        = f(\tau, x(\tau)),
        \quad \tau \in [t, T], \quad x(\tau) \in \mathbb{R},
    \end{equation}
    under the initial condition
    \begin{equation} \label{initial_condition_t}
        x(\tau)
        = w(\tau),
        \quad \tau \in [0, t],
    \end{equation}
    where $T > 0$ is a fixed number, $f \colon [0, T] \times \mathbb{R} \to \mathbb{R}$ is a given function, $(^C D_{0 +}^\alpha x)(\tau)$ denotes the Caputo derivative of order $\alpha \in (0, 1)$, and the pair $(t, w(\cdot))$ consisting of a number $t \in [0, T)$ and a function $w \colon [0, t] \to \mathbb{R}$ is initial data.

    We first make some comments on Cauchy problem~\eqref{differential_equation_t}, \eqref{initial_condition_t}.

    \smallskip

    (i)
        In the particular case when $t = 0$, this problem takes the form
        \begin{equation} \label{differential_equation_0}
            (^C D_{0 +}^\alpha x)(\tau)
            = f(\tau, x(\tau)),
            \quad \tau \in [0, T], \quad x(\tau) \in \mathbb{R},
        \end{equation}
        subject to
        \begin{equation} \label{initial_condition_0}
            x(0)
            = x_0,
        \end{equation}
        where we denote $x_0 \triangleq w(0)$.
        Cauchy problems of the form~\eqref{differential_equation_0}, \eqref{initial_condition_0} are widely considered in the literature (see, e.g.,~\cite{Kilbas_Srivastava_Trujillo_2006,Diethelm_2010}).
        Note that the choice of the initial point $t = 0$ in condition~\eqref{initial_condition_0} reflects the fact that this point is the inferior limit of the Caputo derivative $(^C D^\alpha_{0 +} x)(\tau)$.
        Thus, problem~\eqref{differential_equation_t}, \eqref{initial_condition_t} is a generalization of problem~\eqref{differential_equation_0}, \eqref{initial_condition_0} to the case of initial conditions given at an arbitrary initial point $t \in [0, T)$.

    \smallskip

    (ii)
        If $t \in (0, T)$, then, for a correct statement of a Cauchy problem for differential equation~\eqref{differential_equation_t}, it is not enough to specify only the value $x(t)$, but it is required to specify all the values $x(\tau)$ for $\tau \in [0, t]$, as it is done by initial condition~\eqref{initial_condition_t}.
        The reason for this is that the Caputo fractional derivative has the hereditary nature, which means that the value $(^C D^\alpha_{0 +} x)(\tau)$ depends on all the values $x(\xi)$ for $\xi \in [0, \tau]$.
        This circumstance also indicates a certain similarity between problem~\eqref{differential_equation_t}, \eqref{initial_condition_t} and Cauchy problems typically formulated for functional differential equations.

    \smallskip

    (iii)
        In spite of the fact that the results presented in the paper are obtained for scalar differential equation~\eqref{differential_equation_t}, they admit a straightforward generalization to the vector case.
        In this regard, it should be noted that, if the Cauchy problem is formulated for differential equation~\eqref{differential_equation_0} subject to the condition $x(t) = w$ for some $(t, w) \in [0, T) \times \mathbb{R}$, then this problem is well-posed under mild assumptions, while its multidimensional counterpart is in general ill-posed (see, e.g.,~\cite{Cong_Tuan_2017}).

    \smallskip

    (iv)
        The need to consider Cauchy problems of the form~\eqref{differential_equation_t}, \eqref{initial_condition_t} arises naturally in the study of optimal control problems and differential games for dynamical systems described by Cauchy problems of the form~\eqref{differential_equation_0}, \eqref{initial_condition_0}.
        Namely, when the dynamic programming approach is applied and methods for constructing optimal feedback control strategies are developed.
        For details, the reader is referred to, e.g.,~\cite{Gomoyunov_2020_SIAM,Gomoyunov_2020_Proc,Gomoyunov_2020_ACS,Gomoyunov_2020_DGA}.

    \smallskip

    A question of existence and uniqueness of a solution $x \colon [0, T] \to \mathbb{R}$ of Cauchy problem~\eqref{differential_equation_t}, \eqref{initial_condition_t}, which is denoted by $x(\cdot) \triangleq x(\cdot \mid t, w(\cdot))$, is addressed in~\cite[Proposition~2]{Gomoyunov_2020_DGA}, where also a semigroup property of solutions is established.
    The case when differential equation~\eqref{differential_equation_t} is linear is investigated separately in~\cite{Gomoyunov_2020_FCAA}.
    Continuity properties of the mapping $(t, w(\cdot)) \mapsto x(\cdot \mid t, w(\cdot))$ are examined in the proofs of~\cite[Theorem~1 and Lemma~1]{Gomoyunov_Lukoyanov_2021} (see also~\cite[Assertion 6]{Gomoyunov_2020_Proc} and~\cite[Sections~7 and~8]{Gomoyunov_2020_SIAM}).
    In~\cite[Theorem~2]{Gomoyunov_2020_DE}, some related results are obtained for fractional differential inclusions.
    The main objective of the present paper is to address a more subtle question of differentiability of this mapping.
    More precisely, we deal with the functional
    \begin{equation} \label{rho_introduction}
        (t, w(\cdot)) \mapsto \rho(t, w(\cdot)) \triangleq x(T \mid t, w(\cdot)),
    \end{equation}
    associating initial data $(t, w(\cdot))$ with the endpoint $x(T \mid t, w(\cdot))$ of the solution of Cauchy problem~\eqref{differential_equation_t}, \eqref{initial_condition_t}.

    Let us now discuss what kind of derivatives of this functional are considered in the paper.
    First of all, note that, if the argument $t$ is fixed, then, for the functional $w(\cdot) \mapsto \rho(t, w(\cdot))$, the classical derivatives (e.g., in Gateaux or Fr\'{e}chet sense) are applicable.
    It should be mentioned here that, in the particular case when $t = 0$, i.e., for Cauchy problem~\eqref{differential_equation_0}, \eqref{initial_condition_0}, differentiability properties of the function $x_0 \mapsto x(T \mid 0, x_0)$ are established in~\cite[Proposition~3.11]{Bergounioux_Bourdin_2020}.
    However, the difficulty is that we can not fix the argument $w(\cdot)$ and vary the argument $t$ independently (especially in the positive direction), because the domain of the function $w(\cdot)$ is the interval $[0, t]$.
    A reasonable way to overcome this is to consider an increment of the functional $\rho$ along appropriate extensions of the function $w(\cdot)$.
    In this regard, the framework of coinvariant ($ci$-) derivatives (see, e.g.,~\cite{Kim_1985,Lukoyanov_2000_JAMM}) can be used.
    The $ci$-derivatives have proved to be a convenient and powerful tool for investigating various problems of the theory of functional differential equations.
    For example, when stability issues are studied by the Lyapunov's second method (see, e.g.,~\cite{Kim_1985,Kim_1999}) and the value functionals of optimal control problems and differential games are characterized via the corresponding path-dependent Hamilton--Jacobi equations (see, e.g.,~\cite{Lukoyanov_2000_JAMM,Lukoyanov_2010,Lukoyanov_Gomoyunov_Plaksin_2017_Doklady} and also~\cite{Gomoyunov_Lukoyanov_Plaksin_2021} and the references therein).
    In particular, $ci$-differentiability properties of solutions of functional differential equations with respect to initial data are established in \cite[Theorem~4.3.1]{Kim_1999}.
    The $ci$-derivatives are closely related (see, e.g.,~\cite[Section~5.2]{Gomoyunov_Lukoyanov_Plaksin_2021}) to vertical and horizontal derivatives (see, e.g.,~\cite{Dupire_2009}), which are also widely used in the theory of path-dependent Hamilton--Jacobi equations, especially of those that are associated to optimal control problems and differential games for dynamical systems described by stochastic functional differential equations (see, e.g.,~\cite{Tang_Zhang_2015,Saporito_2019,Zhou_2021} and also~\cite{Ekren_Touzi_Zhang_2016} and the references therein).
    In the present paper, we follow the approach to the notion of $ci$-derivatives developed in, e.g.,~\cite{Lukoyanov_2000_JAMM,Lukoyanov_2010,Lukoyanov_Gomoyunov_Plaksin_2017_Doklady,Gomoyunov_Lukoyanov_Plaksin_2021}, since admissible initial functions $w(\cdot)$ in Cauchy problem~\eqref{differential_equation_t}, \eqref{initial_condition_t} are assumed to be at least continuous.
    Let us emphasize that the definition of the vertical derivative requires dealing with functions $w(\cdot)$ with jump discontinuities.
    In turn, this leads to the need for an appropriate modification of the notion of a solution of Cauchy problem~\eqref{differential_equation_t}, \eqref{initial_condition_t}, which is beyond the scope of this paper.
    Finally, note that, stimulated by the results~\cite{Gomoyunov_2020_SIAM} on Hamilton--Jacobi--Bellman equations associated to optimal control problems for dynamical systems described by fractional differential equations with the Caputo derivatives, we consider so-called fractional $ci$-derivatives of order $\alpha$, which are a certain generalization of the usual notion of $ci$-derivatives.

    The main contribution of the paper is the proof of the fact that the functional $\rho$ from~\eqref{rho_introduction} is $ci$-differentiable of order $\alpha$.
    Furthermore, formulas for calculation of the corresponding derivatives are provided.

    The paper is organized as follows.
    In Section~\ref{section_Statement_of_the_problem}, we describe a space of admissible initial data $(t, w(\cdot))$, specify assumptions on the function $f$, give a precise definition of the functional $\rho$, and recall the notion of $ci$-differentiability of order $\alpha$.
    In Section~\ref{section_main_result}, we formulate the main result of the paper.
    Sections~\ref{section_change_of_variables}--\ref{section_proof} are devoted to a detailed proof of this result.
    In Section~\ref{section_conclusion}, we make some concluding remarks.

\section{Statement of the problem}
\label{section_Statement_of_the_problem}

\setcounter{section}{2}
\setcounter{equation}{0}\setcounter{theorem}{0}

    For every $t \in [0, T]$, following, e.g.,~\cite[Definition~2.3]{Samko_Kilbas_Marichev_1993}, we consider the space $\AC^\alpha[0, t]$ of all functions $x \colon [0, t] \to \mathbb{R}$ each of which can be represented in the form
    \begin{equation} \label{x_f}
        x(\tau)
        = x(0) + \frac{1}{\Gamma(\alpha)} \int_{0}^{\tau} \frac{\ell(\xi)}{(\tau - \xi)^{1 - \alpha}} \, \rd \xi,
        \quad \tau \in [0, t],
    \end{equation}
    for a (Lebesgue) measurable and essentially bounded function $\ell \colon [0, T] \to \mathbb{R}$.\linebreak
    Note that, in the right-hand side of equality~\eqref{x_f}, the second term is the Riemann--Liouville fractional integral of order $\alpha$ of the function $\ell(\cdot)$ (see, e.g.,~\cite[Definition~2.1]{Samko_Kilbas_Marichev_1993}) and $\Gamma$ is the gamma-function.
    According to, e.g.,~\cite[Remark~3.3 and Theorem~2.4]{Samko_Kilbas_Marichev_1993}, every function $x(\cdot) \in \AC^\alpha[0, t]$ is continuous and has at almost every (a.e.) $\tau \in [0, t]$ a Caputo fractional derivative of order $\alpha$, which is defined by (see, e.g.,~\cite[Section~2.4]{Kilbas_Srivastava_Trujillo_2006} and~\cite[Chapter~3]{Diethelm_2010})
    \begin{equation*}
        (^C D^\alpha_{0 +} x)(\tau)
        \triangleq \frac{1}{\Gamma(1 - \alpha)} \frac{\rd}{\rd \tau} \int_{0}^{\tau} \frac{x(\xi) - x(0)}{(\tau - \xi)^\alpha} \, \rd \xi.
    \end{equation*}
    Moreover, if representation~\eqref{x_f} is valid for some measurable and essentially bounded function $\ell(\cdot)$, then $(^C D^\alpha_{0 +} x)(\tau) = \ell(\tau)$ for a.e. $\tau \in [0, t]$.

    We assume that the function $f \colon [0, T] \times \mathbb{R} \to \mathbb{R}$ from the right-hand side of differential equation~\eqref{differential_equation_t} is continuously differentiable and satisfies the following sublinear growth condition:
    there exists $\gamma \geq 0$ such that
    \begin{equation} \label{sublinear_growth}
        |f(\tau, x)|
        \leq \gamma (1 + |x|),
        \quad \tau \in [0, T], \quad x \in \mathbb{R}.
    \end{equation}
    Following, e.g.,~\cite[Section~3.2]{Gomoyunov_2020_DGA} and~\cite[Section~5]{Gomoyunov_2020_SIAM}, as a space of admissible initial data $(t, w(\cdot))$ in Cauchy problem~\eqref{differential_equation_t}, \eqref{initial_condition_t}, we consider the set
    \begin{equation*}
        G
        \triangleq \bigcup_{t \in [0, T)} \bigl( \{ t \} \times \AC^\alpha[0, t] \bigr).
    \end{equation*}
    For every $(t, w(\cdot)) \in G$, by a solution of Cauchy problem~\eqref{differential_equation_t}, \eqref{initial_condition_t}, we mean a function $x(\cdot) \in \AC^\alpha [0, T]$ that meets initial condition~\eqref{initial_condition_t} and satisfies differential equation~\eqref{differential_equation_t} for a.e. $\tau \in [t, T]$.
    By~\cite[Proposition~2]{Gomoyunov_2020_DGA}, such a solution exists and is unique, and we denote it by $x(\cdot) \triangleq x(\cdot \mid t, w(\cdot))$.

    Now, we may give the precise definition of the functional $\rho$ from~\eqref{rho_introduction}:
    \begin{equation} \label{rho}
        \rho(t, w(\cdot))
        \triangleq x(T \mid t, w(\cdot)),
        \quad (t, w(\cdot)) \in G.
    \end{equation}

    Let us recall a notion of $ci$-differentiability of order $\alpha$ of a functional $\varphi \colon G \to \mathbb{R}$, introduced in~\cite[Section~9]{Gomoyunov_2020_SIAM}.
    Take $(t, w(\cdot)) \in G$ and define the set of admissible extensions of the function $w(\cdot)$ on $[0, T]$ by
    \begin{equation} \label{Lambda_t_w}
        \Lambda(t, w(\cdot))
        \triangleq \bigl\{ \lambda(\cdot) \in \AC^\alpha[0, T] \colon
        \, \lambda_t(\cdot) = w(\cdot) \bigr\}.
    \end{equation}
    Here and below, we denote by $\lambda_t(\cdot) \in \AC^\alpha[0, t]$ the restriction of the function $\lambda(\cdot)$ to $[0, t]$, i.e., $\lambda_t(\xi) \triangleq \lambda(\xi)$ for all $\xi \in [0, t]$.
    Then, the functional $\varphi$ is said to be $ci$-differentiable of order $\alpha$ at the point $(t, w(\cdot))$ if there are $\partial_t^\alpha \varphi (t, w(\cdot))$, $\nabla^\alpha \varphi (t, w(\cdot)) \in \mathbb{R}$ such that, for every $\lambda(\cdot) \in \Lambda(t, w(\cdot))$ and every $\tau \in (t, T)$, the relation below holds:
    \begin{align}
        & \varphi(\tau, \lambda_\tau(\cdot)) - \varphi(t, w(\cdot))
        \nonumber \\
        & \ = \partial_t^\alpha \varphi(t, w(\cdot)) (\tau - t)
        + \nabla^\alpha \varphi(t, w(\cdot)) \int_{t}^{\tau} (^C D^\alpha_{0 +} \lambda)(\xi) \, \rd \xi
        + o(\tau - t),
        \label{ci-differentiability}
    \end{align}
    where the function $o \colon (0, \infty) \to \mathbb{R}$, which may depend on $t$ and $\lambda(\cdot)$, satisfies the condition $\frac{o(\delta)}{\delta} \to 0$ as $\delta \to 0^+$.
    In this case, the values $\partial_t^\alpha \varphi (t, w(\cdot))$ and $\nabla^\alpha \varphi(t, w(\cdot))$ are called $ci$-derivatives of order $\alpha$ of $\varphi$ at $(t, w(\cdot))$.

\section{Main result}
\label{section_main_result}

\setcounter{section}{3}
\setcounter{equation}{0}\setcounter{theorem}{0}

    Fix $(t, w(\cdot)) \in G$.
    For every $\tau \in (t, T]$, denote
    \begin{equation} \label{overline_p_q}
        \overline{p}(\tau)
        \triangleq - \frac{1 - \alpha}{\Gamma(\alpha)}  \biggl( \frac{T - \tau}{T - t} \biggr)
        \int_{0}^{t} \frac{(^C D^\alpha_{0 +} w)(\xi)}{(\tau - \xi)^{2 - \alpha}} \, \rd \xi,
        \quad \overline{q}(\tau)
        \triangleq \frac{1}{\Gamma(\alpha) (\tau - t)^{1 - \alpha}}.
    \end{equation}
    Take the solution $x(\cdot) \triangleq x(\cdot \mid t, w(\cdot))$ of Cauchy problem~\eqref{differential_equation_t}, \eqref{initial_condition_t}.
    Put
    \begin{equation} \label{a_b}
        a(\tau)
        \triangleq \frac{\partial f}{\partial x} (\tau, x(\tau)),
        \quad b(\tau)
        \triangleq \biggl( \frac{T - \tau}{T - t} \biggr) \frac{\partial f}{\partial \tau} (\tau, x(\tau))
        - \frac{\alpha f(\tau, x(\tau))}{T - t}
    \end{equation}
    for all $\tau \in [t, T]$, where $\frac{\partial f}{\partial \tau}$ and $\frac{\partial f}{\partial x}$ are the partial derivatives of the function $f$.
    If $\tau = 0$ or $\tau = T$, then the partial derivative $\frac{\partial f}{\partial \tau}$ is understood as the corresponding one-sided derivative.

    Consider the linear weakly-singular Volterra integral equations
    \begin{equation} \label{integral_equation_p}
        p(\tau)
        = \overline{p}(\tau)
        + \frac{1}{\Gamma(\alpha)} \int_{t}^{\tau} \frac{a(\xi) p(\xi) + b(\xi)}{(\tau - \xi)^{1 - \alpha}} \, \rd \xi,
        \quad \tau \in (t, T],
    \end{equation}
    and
    \begin{equation} \label{integral_equation_q}
        q(\tau)
        = \overline{q}(\tau)
        + \frac{1}{\Gamma(\alpha)} \int_{t}^{\tau} \frac{a(\xi) q(\xi)}{(\tau - \xi)^{1 - \alpha}} \, \rd \xi,
        \quad \tau \in (t, T].
    \end{equation}
    By a solution of integral equation~\eqref{integral_equation_p} (respectively, integral equation~\eqref{integral_equation_q}), we mean a function $p(\cdot)$ (respectively, a function $q(\cdot)$) from $\C^{1 - \alpha}(t, T]$ that satisfies this equation for every $\tau \in (t, T]$.
    Here, we denote by $\C^{1 - \alpha} (t, T]$ the space of all continuous functions $p \colon (t, T] \to \mathbb{R}$ such that the function $(t, T] \ni \tau \mapsto (\tau - t)^{1 - \alpha} p(\tau)$ is bounded.

    \begin{theo} \label{theorem}
        The functional $\rho$ defined by~\eqref{rho} is $ci$-differentiable of order $\alpha$ at every point $(t, w(\cdot)) \in G$, and the corresponding derivatives are given by
        \begin{equation*}
            \partial_t^\alpha \rho(t, w(\cdot))
            = p(T),
            \quad \nabla^\alpha \rho(t, w(\cdot))
            = q(T),
        \end{equation*}
        where $p(\cdot) \triangleq p(\cdot \mid t, w(\cdot))$ and $q(\cdot) \triangleq q(\cdot \mid t, w(\cdot))$ are the unique solutions of integral equations~\eqref{integral_equation_p} and~\eqref{integral_equation_q}, respectively.
    \end{theo}

    The rest of the paper is devoted to the proof of Theorem~\ref{theorem}.

    As a first step, in the next section, we move from Cauchy problem~\eqref{differential_equation_t}, \eqref{initial_condition_t} to the corresponding integral equation and make a change of variables in order to obtain a unified interval $[0, 1]$ instead of the interval $[t, T]$.

\section{Equivalent integral equation and change of variables}
\label{section_change_of_variables}

\setcounter{section}{4}
\setcounter{equation}{0}\setcounter{theorem}{0}

    Fix $(t, w(\cdot)) \in G$.
    According to, e.g.,~\cite[Proposition~2]{Gomoyunov_2020_DGA}, the solution $x(\cdot) \triangleq x(\cdot \mid t, w(\cdot))$ of Cauchy problem~\eqref{differential_equation_t}, \eqref{initial_condition_t} is the unique function from $\C[0, T]$,  the space of all continuous functions from $[0, T]$ to $\mathbb{R}$, that satisfies initial condition~\eqref{initial_condition_t} and the integral equation
    \begin{equation} \label{integral_equation}
        x(\tau)
        = \overline{x}(\tau)
        + \frac{1}{\Gamma(\alpha)} \int_{t}^{\tau} \frac{f(\xi, x(\xi))}{(\tau - \xi)^{1 - \alpha}} \, \rd \xi,
        \quad \tau \in [t, T],
    \end{equation}
    where the function $\overline{x}(\cdot) \triangleq \overline{x}(\cdot \mid t, w(\cdot))$ is defined by
    \begin{equation} \label{overline_x}
        \overline{x}(\tau)
        \triangleq
        \begin{cases}
            w(\tau), & \mbox{if } \tau \in [0, t], \\
            \displaystyle
            w(0) + \frac{1}{\Gamma(\alpha)} \int_{0}^{t} \frac{(^C D^\alpha_{0 +} w)(\xi)}{(\tau - \xi)^{1 - \alpha}} \, \rd \xi, & \mbox{if } \tau \in (t, T].
          \end{cases}
    \end{equation}
    Note that $\overline{x}(\cdot) \in \AC^\alpha[0, T]$ and $(^C D^\alpha_{0 +} \overline{x})(\tau) = 0$ for all $\tau \in (t, T)$.

    Making the change of variables $\vartheta \triangleq \frac{\tau - t}{T - t}$ in integral equation~\eqref{integral_equation}, we come to the auxiliary integral equation
    \begin{equation} \label{y_integral_equation}
        \mathbf{x}(\vartheta)
        = \overline{\mathbf{x}}(\vartheta)
        + \frac{(T - t)^\alpha}{\Gamma(\alpha)} \int_{0}^{\vartheta}
        \frac{f(t + \zeta (T - t), \mathbf{x}(\zeta))}{(\vartheta - \zeta)^{1 - \alpha}} \, \rd \zeta,
        \quad \vartheta \in [0, 1],
    \end{equation}
    where the function $\overline{\mathbf{x}}(\cdot) \triangleq \overline{\mathbf{x}}(\cdot \mid t, w(\cdot))$ is given by
    \begin{equation} \label{overline_mathbf_x}
        \overline{\mathbf{x}}(\vartheta)
        \triangleq \overline{x}(t + \vartheta (T - t) \mid t, w(\cdot)),
        \quad \vartheta \in [0, 1].
    \end{equation}
    By a solution of integral equation~\eqref{y_integral_equation}, we mean a function $\mathbf{x}(\cdot) \in \C[0, 1]$ satisfying this equation.
    Here and below, we use bold letters to denote functions of the new variable $\vartheta$.

    \begin{lem} \label{lemma_mathbf_x_x}
        There is a unique solution $\mathbf{x}(\cdot) \triangleq \mathbf{x}(\cdot \mid t, w(\cdot))$ of integral equation~\eqref{y_integral_equation}.
        Moreover, the following equality holds:
        \begin{equation} \label{lemma_mathbf_x_x_main}
            \mathbf{x}(\vartheta)
            = x(t + \vartheta (T - t)),
            \quad \vartheta \in [0, 1],
        \end{equation}
        where $x(\cdot) \triangleq x(\cdot \mid t, w(\cdot))$ is the solution of Cauchy problem~\eqref{differential_equation_t}, \eqref{initial_condition_t}.
    \end{lem}
    \proof
        On the one hand, let us define a function $\mathbf{x} \colon [0, 1] \to \mathbb{R}$ by $x(\cdot)$ according to~\eqref{lemma_mathbf_x_x_main}.
        Since the function $x(\cdot)$ is continuous and satisfies integral equation~\eqref{integral_equation}, we have $\mathbf{x}(\cdot) \in \C[0, 1]$ and, for all $\vartheta \in [0, 1]$,
        \begin{align*}
            & \mathbf{x}(\vartheta)
            = \overline{x}(t + \vartheta (T - t))
            + \frac{1}{\Gamma(\alpha)} \int_{t}^{t + \vartheta (T - t)} \frac{f(\xi, x(\xi))}{(t + \vartheta (T - t) - \xi)^{1 - \alpha}} \, \rd \xi \\
            & \hphantom{\mathbf{x}(\vartheta)}
            = \overline{\mathbf{x}}(\vartheta)
            + \frac{(T - t)^\alpha}{\Gamma(\alpha)} \int_{0}^{\vartheta} \frac{f(t + \zeta (T - t), \mathbf{x}(\zeta))}
            {(\vartheta - \zeta)^{1 - \alpha}} \, \rd \zeta.
        \end{align*}
        Hence, the function $\mathbf{x}(\cdot)$ is a solution of integral equation~\eqref{y_integral_equation}.

        On the other hand, let us take a solution $\mathbf{x}(\cdot)$ of integral equation~\eqref{y_integral_equation} and consider the function
        \begin{equation} \label{proposition_y_x_proof_widetilde_x}
            \widetilde{x}(\tau)
            \triangleq
            \begin{cases}
                w(\tau), & \mbox{if } \tau \in [0, t), \\
                \displaystyle
                \mathbf{x} \biggl( \frac{\tau - t}{T - t} \biggr), & \mbox{if } \tau \in [t, T].
            \end{cases}
        \end{equation}
        Since $\mathbf{x}(0) = \overline{\mathbf{x}}(0) = \overline{x}(t) = w(t)$ and the functions $w(\cdot)$ and $\mathbf{x}(\cdot)$ are continuous, we obtain that the function $\widetilde{x}(\cdot)$ is continuous and meets initial condition~\eqref{initial_condition_t}.
        Furthermore, for every $\tau \in [t, T]$, we derive
        \begin{align*}
            & \widetilde{x}(\tau)
            = \overline{\mathbf{x}}\biggl( \frac{\tau - t}{T - t} \biggr)
            + \frac{(T - t)^\alpha}{\Gamma(\alpha)} \int_{0}^{\frac{\tau - t}{T - t}}
            \frac{f(t + \zeta (T - t), \mathbf{x}(\zeta))}{(\frac{\tau - t}{T - t} - \zeta)^{1 - \alpha}} \, \rd \zeta \\
            & \hphantom{\widetilde{x}(\tau)}
            = \overline{x}(\tau)
            + \frac{1}{\Gamma(\alpha)} \int_{t}^{\tau}
            \frac{f(\xi, \widetilde{x}(\xi))}{(\tau - \xi)^{1 - \alpha}} \, \rd \xi,
        \end{align*}
        which means that the function $\widetilde{x}(\cdot)$ satisfies integral equation~\eqref{integral_equation}.
        Thus, we have $\widetilde{x}(\cdot) = x(\cdot)$, and, therefore, due to~\eqref{proposition_y_x_proof_widetilde_x}, the result follows.
    \proofend

    In particular, Lemma~\ref{lemma_mathbf_x_x} implies that the functional $\rho$ from~\eqref{rho} can be represented in the following way:
    \begin{equation} \label{varphi_y}
        \rho(t, w(\cdot))
        = \mathbf{x}(1 \mid t, w(\cdot)),
        \quad (t, w(\cdot)) \in G.
    \end{equation}

    To examine properties of the solution $\mathbf{x}(\cdot \mid t, w(\cdot))$ of auxiliary integral equation~\eqref{y_integral_equation}, we start by investigating the free term $\overline{\mathbf{x}}(\cdot \mid t, w(\cdot))$ of this equation separately.

\section{Properties of free term of auxiliary integral equation}
\label{section_properties_free_term}

\setcounter{section}{5}
\setcounter{equation}{0}\setcounter{theorem}{0}

    Let $(t, w(\cdot)) \in G$ and $\ell \in \mathbb{R}$ be fixed.
    Consider the extension $\lambda^{(\ell)}(\cdot) \triangleq \lambda^{(\ell)}(\cdot \mid t, w(\cdot)) \in \Lambda(t, w(\cdot))$ (see~\eqref{Lambda_t_w}) such that $(^C D^\alpha_{0 +} \lambda^{(\ell)})(\tau) = \ell$ for all $\tau \in (t, T)$.
    Note that, in explicit form, we have
    \begin{equation}\label{x^f}
        \lambda^{(\ell)}(\tau)
        = \begin{cases}
            w(\tau), & \mbox{if } \tau \in [0, t], \\
            \displaystyle
            w(0) + \frac{1}{\Gamma(\alpha)} \int_{0}^{t} \frac{(^C D^\alpha_{0 +} w)(\xi)}{(\tau - \xi)^{1 - \alpha}} \, \rd \xi
            + \frac{\ell (\tau - t)^\alpha}{\Gamma(\alpha)}, & \mbox{if } \tau \in (t, T].
          \end{cases}
    \end{equation}
    For every $\tau \in [t, T)$, we take the function $\overline{\mathbf{x}}(\cdot \mid \tau, \lambda^{(\ell)}_\tau(\cdot))$ corresponding to the pair $(\tau, \lambda^{(\ell)}_\tau(\cdot))$ (see~\eqref{overline_mathbf_x}) and put $\overline{\mathbf{x}}(\cdot \mid \tau) \triangleq \overline{\mathbf{x}}(\cdot \mid \tau, \lambda^{(\ell)}_\tau(\cdot))$ for brevity.

    \begin{lem} \label{lemma_overline_mathbf_Lipschitz}
        For any $\eta \in (0, T - t)$, there exists $\overline{L} \geq 0$ such that, for every $\tau \in (t, T - \eta]$ and every $\vartheta \in (0, 1]$, the inequality below is fulfilled:
        \begin{equation} \label{lemma_overline_mathbf_Lipschitz_main}
            \biggl| \frac{\overline{\mathbf{x}}(\vartheta \mid \tau) - \overline{\mathbf{x}}(\vartheta \mid t)}{\tau - t} \biggr|
            \leq \frac{\overline{L}}{\vartheta^{1 - \alpha}}.
        \end{equation}
    \end{lem}
    \proof
        Choose $M \geq 0$ such that
        \begin{equation} \label{M}
            |(^C D^\alpha_{0 +} w)(\xi)|
            \leq M
            \text{ for a.e. } \xi \in [0, t]
        \end{equation}
        and, for a given $\eta \in (0, T - t)$, define
        \begin{equation*}
            \overline{L} \triangleq \frac{M + |\ell|}{\Gamma(\alpha) \eta^{1 - \alpha}}.
        \end{equation*}

        Fix $\tau \in (t, T - \eta]$ and $\vartheta \in (0, 1]$.
        According to~\eqref{overline_x}, \eqref{overline_mathbf_x}, and~\eqref{x^f}, we have
        \begin{align*}
            & |\overline{\mathbf{x}}(\vartheta \mid \tau) - \overline{\mathbf{x}}(\vartheta \mid t)| \\
            & \quad \leq \frac{M}{\Gamma(\alpha)}
            \int_{0}^{t} \biggl( \frac{1}{(t + \vartheta (T - t) - \xi)^{1 - \alpha}}
            - \frac{1}{(\tau + \vartheta (T - \tau) - \xi)^{1 - \alpha}} \biggr) \, \rd \xi \\
            & \quad + \frac{|\ell|}{\Gamma(\alpha)} \int_{t}^{\tau} \frac{\rd \xi}{(\tau + \vartheta (T - \tau) - \xi)^{1 - \alpha}}.
        \end{align*}
        Then, in view of the estimates
        \begin{align*}
            & \int_{0}^{t} \biggl( \frac{1}{(t + \vartheta (T - t) - \xi)^{1 - \alpha}}
            - \frac{1}{(\tau + \vartheta (T - \tau) - \xi)^{1 - \alpha}} \biggr) \, \rd \xi \\
            & \quad \leq (1 - \alpha) (1 - \vartheta) (\tau - t) \int_{0}^{t} \frac{\rd \xi}{(t + \vartheta (T - t) - \xi)^{2 - \alpha}} \\
            & \quad \leq \frac{(1 - \vartheta)(\tau - t)}{\vartheta^{1 - \alpha} (T - t)^{1 - \alpha}}
            \leq \frac{\tau - t}{\vartheta^{1 - \alpha} \eta^{1 - \alpha}}
        \end{align*}
        and
        \begin{equation*}
            \int_{t}^{\tau} \frac{\rd \xi}{(\tau + \vartheta (T - \tau) - \xi)^{1 - \alpha}} \\
            \leq \frac{\tau - t}{\vartheta^{1 - \alpha} (T - \tau)^{1 - \alpha}}
            \leq \frac{\tau - t}{\vartheta^{1 - \alpha} \eta^{1 - \alpha}},
        \end{equation*}
        we obtain that inequality~\eqref{lemma_overline_mathbf_Lipschitz_main} is fulfilled for the specified number $\overline{L}$.
    \proofend

    By the functions $\overline{p}(\cdot)$ and $\overline{q}(\cdot)$ from~\eqref{overline_p_q}, define
    \begin{equation} \label{overline_mathbf_p_q}
        \overline{\mathbf{p}}(\vartheta)
        \triangleq \overline{p}(t + \vartheta (T - t)),
        \quad \overline{\mathbf{q}}(\vartheta)
        \triangleq \overline{q}(t + \vartheta (T - t)),
        \quad \vartheta \in (0, 1].
    \end{equation}

    Denote by $\C^{1 - \alpha}(0, 1]$ the space of all continuous functions $\mathbf{p} \colon (0, 1] \to \mathbb{R}$ such that the function $(0, 1] \ni \vartheta \mapsto \vartheta^{1 - \alpha} \mathbf{p}(\vartheta)$ is bounded.

    \begin{lem} \label{lemma_inclusions}
        The inclusions $\overline{\mathbf{p}}(\cdot)$, $\overline{\mathbf{q}}(\cdot) \in \C^{1 - \alpha} (0, 1]$ are valid.
    \end{lem}
    \proof
        Note that the inclusion $\overline{\mathbf{q}}(\cdot) \in \C^{1 - \alpha} (0, 1]$ obviously holds.
        Further, by the theorem on continuity of integrals with respect to a parameter (see, e.g.,~\cite[Corollary 2.8.7, item (i)]{Bogachev_2007}), we find that the function $\overline{\mathbf{p}}(\cdot)$ is continuous.
        In addition, taking the number $M$ from~\eqref{M}, we get
        \begin{align*}
            & |\overline{\mathbf{p}}(\vartheta)|
            \leq \frac{(1 - \alpha) (1 - \vartheta) M}{\Gamma(\alpha)}
            \int_{0}^{t} \frac{\rd \xi}{(t + \vartheta (T - t) - \xi)^{2 - \alpha}} \\
            & \hphantom{|\overline{\mathbf{p}}(\vartheta)|}
            \leq \frac{M}{\Gamma(\alpha) \vartheta^{1 - \alpha} (T - t)^{1 - \alpha}}
        \end{align*}
        for all $\vartheta \in (0, 1]$,
        which implies that the function $(0, 1] \ni \vartheta \mapsto \vartheta^{1 - \alpha} \overline{\mathbf{p}}(\vartheta)$ is bounded.
        Hence, $\overline{\mathbf{p}}(\cdot) \in \C^{1 - \alpha} (0, 1]$, and the proof is complete.
    \proofend

    It should be emphasized that the function $(0, 1] \ni \vartheta \mapsto \vartheta^{1 - \alpha} \overline{\mathbf{p}}(\vartheta)$ may not have a limit as $\vartheta \to 0^+$.
    An example is provided in Appendix~\ref{appendix}.

    \begin{lem} \label{lemma_overline_mathbf_x_derivatives}
        For every $\vartheta \in (0, 1]$, the following limit relations hold:
        \begin{equation} \label{proposition_a_pointwise_limit}
            \lim_{\tau \to t^+}
            \frac{\overline{\mathbf{x}}(\vartheta \mid \tau) - \overline{\mathbf{x}}(\vartheta \mid t)}{\tau - t}
            = \overline{\mathbf{p}}(\vartheta) + \overline{\mathbf{q}}(\vartheta) \ell
        \end{equation}
        and
        \begin{equation} \label{proposition_a_integral_limit}
            \lim_{\tau \to t^+} \int_{0}^{\vartheta}
            \biggl| \frac{\overline{\mathbf{x}}(\zeta \mid \tau) - \overline{\mathbf{x}}(\zeta \mid t)}{\tau - t}
            - \overline{\mathbf{p}}(\zeta) - \overline{\mathbf{q}}(\zeta) \ell \biggr|
            \frac{\rd \zeta}{(\vartheta - \zeta)^{1 - \alpha}}
            = 0.
        \end{equation}
    \end{lem}
    \proof
        Note that relation~\eqref{proposition_a_pointwise_limit} can be verified by direct calculation and based on the theorem on differentiability of integrals with respect to a parameter (see, e.g.,~\cite[Corollary~2.8.7, item (ii)]{Bogachev_2007}).
        Nevertheless, since we need to prove also relation~\eqref{proposition_a_integral_limit}, it is still required to estimate the value
        \begin{equation} \label{overline_mathbf_z}
            \overline{\mathbf{z}}(\vartheta \mid \tau)
            \triangleq \biggl| \frac{\overline{\mathbf{x}}(\vartheta \mid \tau) - \overline{\mathbf{x}}(\vartheta \mid t)}{\tau - t}
            - \overline{\mathbf{p}}(\vartheta) - \overline{\mathbf{q}}(\vartheta) \ell \biggr|
        \end{equation}
        for every $\vartheta \in (0, 1]$ and every $\tau \in (t, T)$.

        In accordance with the notation introduced in~\eqref{overline_x}, \eqref{overline_mathbf_x}, \eqref{x^f} and also in~\eqref{overline_p_q}, \eqref{overline_mathbf_p_q}, taking the number $M$ from~\eqref{M}, we derive
        \begin{align*}
            & \overline{\mathbf{z}}(\vartheta \mid \tau)(\tau - t)
            \leq \frac{M (1 - \alpha) (1 - \vartheta) (\tau - t)}{\Gamma(\alpha)} \int_{0}^{t} \frac{\rd \xi}{(t + \vartheta (T - t) - \xi)^{2 - \alpha}} \\
            & \quad + \frac{M}{\Gamma(\alpha)} \int_{0}^{t} \biggl( \frac{1}{(\tau + \vartheta (T - \tau) - \xi)^{1 - \alpha}}
            - \frac{1}{(t + \vartheta (T - t) - \xi)^{1 - \alpha}} \biggr) \, \rd \xi \\
            & \quad + \frac{|\ell|}{\Gamma(\alpha)} \biggl| \int_{t}^{\tau} \frac{\rd \xi}{(\tau + \vartheta (T - \tau) - \xi)^{1 - \alpha}}
            - \frac{\tau - t}{\vartheta^{1 -\alpha} (T - t)^{1 - \alpha}} \biggr|.
        \end{align*}
        We have
        \begin{align*}
            & \int_{0}^{t} \biggl( \frac{1}{(\tau + \vartheta (T - \tau) - \xi)^{1 - \alpha}}
            - \frac{1}{(t + \vartheta (T - t) - \xi)^{1 - \alpha}} \biggr) \, \rd \xi \\
            & \quad \leq - (1 - \alpha) (1 - \vartheta) (\tau - t) \int_{0}^{t} \frac{\rd \xi}{(\tau + \vartheta (T - \tau) - \xi)^{2 - \alpha}}
        \end{align*}
        and
        \begin{align*}
            & \int_{0}^{t} \frac{\rd \xi}{(t + \vartheta (T - t) - \xi)^{2 - \alpha}}
            - \int_{0}^{t} \frac{\rd \xi}{(\tau + \vartheta (T - \tau) - \xi)^{2 - \alpha}} \\
            & \quad \leq \frac{1}{1 - \alpha} \biggl( \frac{1}{\vartheta^{1 - \alpha} (T - t)^{1 - \alpha}}
            - \frac{1}{(\tau + \vartheta (T - \tau) - t)^{1 - \alpha}} \biggr).
        \end{align*}
        In addition, we obtain
        \begin{align*}
            & \biggl| \int_{t}^{\tau} \frac{\rd \xi}{(\tau + \vartheta (T - \tau) - \xi)^{1 - \alpha}}
            - \frac{\tau - t}{\vartheta^{1 - \alpha} (T - t)^{1 - \alpha}} \biggr| \\
            & \quad \leq \biggl( \frac{\tau - t}{\vartheta^{1 - \alpha} (T - \tau)^{1 - \alpha}}
            - \int_{t}^{\tau} \frac{\rd \xi}{(\tau + \vartheta (T - \tau) - \xi)^{1 - \alpha}} \biggr) \\
            & \quad + \biggl( \frac{\tau - t}{\vartheta^{1 - \alpha} (T - \tau)^{1 - \alpha}}
            - \frac{\tau - t}{\vartheta^{1 - \alpha} (T - t)^{1 - \alpha}} \biggr) \\
            & \quad \leq 2 (\tau - t) \biggl( \frac{1}{\vartheta^{1 - \alpha} (T - \tau)^{1 - \alpha}}
            - \frac{1}{(\tau + \vartheta (T - \tau) - t)^{1 - \alpha}} \biggr).
        \end{align*}
        As a result, we get
        \begin{equation} \label{proposition_a_proof_main_estimate}
            \overline{\mathbf{z}}(\vartheta \mid \tau)
            \leq \frac{M + 2 |\ell|}{\Gamma(\alpha)} \biggl( \frac{1}{\vartheta^{1 - \alpha} (T - \tau)^{1 - \alpha}}
            - \frac{1}{(\tau + \vartheta (T - \tau) - t)^{1 - \alpha}} \biggr).
        \end{equation}
        In particular, noting that the right-hand side of inequality~\eqref{proposition_a_proof_main_estimate} tends to $0$ as $\tau \to t^+$, we find that~\eqref{proposition_a_pointwise_limit} is valid.

        Further, due to~\eqref{proposition_a_proof_main_estimate}, we have
        \begin{align*}
            & \int_{0}^{\vartheta} \frac{\overline{\mathbf{z}}(\zeta \mid \tau)}{(\vartheta - \zeta)^{1 - \alpha}} \, \rd \zeta
            \leq \frac{M + 2 |\ell|}{\Gamma(\alpha) (T - \tau)^{1 - \alpha}} \\
            & \quad \times \biggl( \int_{0}^{\vartheta} \frac{\rd \zeta}{\zeta^{1 - \alpha} (\vartheta - \zeta)^{1 - \alpha}}
            - \int_{0}^{\vartheta} \frac{\rd \zeta}{\bigl(\frac{\tau - t}{T - \tau} + \zeta \bigr)^{1 - \alpha} (\vartheta - \zeta)^{1 - \alpha}} \biggr).
        \end{align*}
        Consequently, relying on the relations
        \begin{equation} \label{formula_beta}
            \int_{0}^{\vartheta} \frac{\rd \zeta}{\zeta^{1 - \alpha} (\vartheta - \zeta)^{1 - \alpha}}
            = B(\alpha, \alpha) \vartheta^{2 \alpha - 1}
        \end{equation}
        and
        \begin{align*}
            & \int_{0}^{\vartheta} \frac{\rd \zeta}{\bigl(\frac{\tau - t}{T - \tau} + \zeta \bigr)^{1 - \alpha} (\vartheta - \zeta)^{1 - \alpha}} \\
            & \quad = \int_{- \frac{\tau - t}{T - \tau}}^{\vartheta}
            \frac{\rd \zeta}{\bigl(\frac{\tau - t}{T - \tau} + \zeta \bigr)^{1 - \alpha} (\vartheta - \zeta)^{1 - \alpha}}
            - \int_{- \frac{\tau - t}{T - \tau}}^{0}
            \frac{\rd \zeta}{\bigl(\frac{\tau - t}{T - \tau} + \zeta \bigr)^{1 - \alpha} (\vartheta - \zeta)^{1 - \alpha}} \\
            & \quad \geq B(\alpha, \alpha) \biggl( \vartheta + \frac{\tau - t}{T - \tau} \biggr)^{2 \alpha - 1}
            - \frac{1}{\vartheta^{1 - \alpha}} \int_{- \frac{\tau - t}{T - \tau}}^{0}
            \frac{\rd \zeta}{\bigl(\frac{\tau - t}{T - \tau} + \zeta \bigr)^{1 - \alpha}} \\
            & \quad = B(\alpha, \alpha) \biggl( \vartheta + \frac{\tau - t}{T - \tau} \biggr)^{2 \alpha - 1}
            - \frac{1}{\alpha \vartheta^{1 - \alpha}} \biggl( \frac{\tau - t}{T - \tau} \biggr)^\alpha,
        \end{align*}
        where $B$ is the beta-function, we derive
        \begin{align*}
            & \int_{0}^{\vartheta}
            \frac{\overline{\mathbf{z}}(\zeta \mid \tau)}{(\vartheta - \zeta)^{1 - \alpha}} \, \rd \zeta
            \leq \frac{M + 2 |\ell|}{\Gamma(\alpha) (T - \tau)^{1 - \alpha}} \\
            & \quad \times \biggl( B(\alpha, \alpha) \vartheta^{2 \alpha - 1}
            - B(\alpha, \alpha) \biggl( \vartheta + \frac{\tau - t}{T - \tau} \biggr)^{2 \alpha - 1}
            + \frac{1}{\alpha \vartheta^{1 - \alpha}} \biggl( \frac{\tau - t}{T - \tau} \biggr)^\alpha \biggr).
        \end{align*}
        Since the right-hand side of this inequality tends to $0$ as $\tau \to t^+$, we conclude that~\eqref{proposition_a_integral_limit} holds.
        The lemma is proved.
    \proofend

    Before proceeding with the study of the solution of auxiliary integral equation~\eqref{y_integral_equation}, we provide some technical results that concern the unique solvability of a class of linear weakly-singular Volterra integral equations in the space $\C^{1 - \alpha}(0, 1]$ and the corresponding Gronwall-type inequalities.

\section{Linear weakly-singular Volterra integral equations}
\label{section_Volterra_equations}

\setcounter{section}{6}
\setcounter{equation}{0}\setcounter{theorem}{0}

    Let $\overline{\mathbf{y}}(\cdot) \in \C^{1 - \alpha}(0, 1]$ and $\mathbf{c}(\cdot)$, $\mathbf{d}(\cdot) \in \C[0, 1]$ be given.
    Consider the integral equation
    \begin{equation} \label{integral_equation_abstract}
        \mathbf{y}(\vartheta)
        = \overline{\mathbf{y}}(\vartheta)
        + \frac{1}{\Gamma(\alpha)} \int_{0}^{\vartheta} \frac{\mathbf{c}(\zeta) \mathbf{y}(\zeta) + \mathbf{d}(\zeta)}{(\vartheta - \zeta)^{1 - \alpha}} \, \rd \zeta,
        \quad \vartheta \in (0, 1].
    \end{equation}
    A solution of integral equation~\eqref{integral_equation_abstract} is defined as a function $\mathbf{y}(\cdot) \in \C^{1 - \alpha} (0, 1]$ that satisfies this equation.
    \begin{prop} \label{proposition_linear_existence_uniqueness_arbitrary}
        Integral equation~\eqref{integral_equation_abstract} has a unique solution.
    \end{prop}
    \proof
        Suppose that $\mathbf{y}(\cdot)$ is a solution of~\eqref{integral_equation_abstract}.
        Define
        \begin{equation} \label{proposition_linear_existence_uniqueness_arbitrary_proof_z_definition}
            \mathbf{s}(\vartheta)
            \triangleq
            \begin{cases}
                0, & \mbox{if } \vartheta = 0, \\
                \vartheta^{1 - \alpha} (\mathbf{y}(\vartheta) - \overline{\mathbf{y}}(\vartheta)), & \mbox{if } \vartheta \in (0, 1].
            \end{cases}
        \end{equation}
        Then, the function $\mathbf{s}(\cdot)$ is continuous on $(0, 1]$, bounded, and, moreover,
        \begin{equation*}
            \mathbf{s}(\vartheta)
            = \frac{\vartheta^{1 - \alpha}}{\Gamma(\alpha)} \int_{0}^{\vartheta} \frac{\mathbf{c}(\zeta) \mathbf{s}(\zeta)}{\zeta^{1 - \alpha} (\vartheta - \zeta)^{1 - \alpha}} \, \rd \zeta
            + \frac{\vartheta^{1 - \alpha}}{\Gamma(\alpha)} \int_{0}^{\vartheta} \frac{\mathbf{c}(\zeta) \overline{\mathbf{y}}(\zeta) + \mathbf{d}(\zeta)}{(\vartheta - \zeta)^{1 - \alpha}} \, \rd \zeta
        \end{equation*}
        for every $\vartheta \in (0, 1]$.
        In addition, since the functions $(0, 1] \ni \vartheta \mapsto \mathbf{c}(\vartheta) \mathbf{s}(\vartheta)$ and $(0, 1] \ni \vartheta \mapsto \vartheta^{1 - \alpha} \mathbf{c}(\vartheta) \overline{\mathbf{y}}(\vartheta) + \vartheta^{1 - \alpha} \mathbf{d}(\vartheta)$ are continuous and bounded, it follows from~\cite[Corollary~2.1]{Gomoyunov_2020_FCAA} that the functions
        \begin{equation*}
            [0, 1] \ni \vartheta \mapsto
            \frac{\vartheta^{1 - \alpha}}{\Gamma(\alpha)} \int_{0}^{\vartheta} \frac{\mathbf{c}(\zeta) \mathbf{s}(\zeta)}{\zeta^{1 - \alpha} (\vartheta - \zeta)^{1 - \alpha}} \, \rd \zeta
        \end{equation*}
        and
        \begin{equation} \label{proof_proposition_linear_existence_uniqueness_arbitrary_widetilde_a}
            [0, 1] \ni \vartheta \mapsto \overline{\mathbf{s}}(\vartheta)
            \triangleq \frac{\vartheta^{1 - \alpha}}{\Gamma(\alpha)} \int_{0}^{\vartheta} \frac{\mathbf{c}(\zeta) \overline{\mathbf{y}}(\zeta) + \mathbf{d}(\zeta)}{(\vartheta - \zeta)^{1 - \alpha}} \, \rd \zeta
        \end{equation}
        are continuous and tend to $0$ as $\vartheta \to 0^+$.
        Hence, we find that the function $\mathbf{s}(\cdot)$ is continuous on $[0, 1]$ and satisfies the integral equation
        \begin{equation} \label{proof_proposition_linear_existence_uniqueness_arbitrary_auxiliaty_integral_equation}
            \mathbf{s}(\vartheta)
            = \overline{\mathbf{s}}(\vartheta)
            + \frac{\vartheta^{1 - \alpha}}{\Gamma(\alpha)}
            \int_{0}^{\vartheta} \frac{\mathbf{c}(\zeta) \mathbf{s}(\zeta)}{\zeta^{1 - \alpha} (\vartheta - \zeta)^{1 - \alpha}} \, \rd \zeta,
            \quad \vartheta \in [0, 1].
        \end{equation}

        Now, let us verify the converse statement.
        Let a function $\mathbf{s}(\cdot) \in \C[0, 1]$ satisfy integral equation~\eqref{proof_proposition_linear_existence_uniqueness_arbitrary_auxiliaty_integral_equation}.
        In accordance with~\eqref{proposition_linear_existence_uniqueness_arbitrary_proof_z_definition}, put
        \begin{equation*}
            \mathbf{y}(\vartheta)
            \triangleq \overline{\mathbf{y}}(\vartheta) + \frac{\mathbf{s}(\vartheta)}{\vartheta^{1 - \alpha}},
            \quad \vartheta \in (0, 1].
        \end{equation*}
        Then, we have $\mathbf{y}(\cdot) \in \C^{1 - \alpha} (0, 1]$, and, taking~\eqref{proof_proposition_linear_existence_uniqueness_arbitrary_widetilde_a} into account, we obtain
        \begin{align*}
            & \mathbf{y}(\vartheta)
            = \overline{\mathbf{y}}(\vartheta)
            + \frac{\overline{\mathbf{s}}(\vartheta)}{\vartheta^{1 - \alpha}}
            + \frac{1}{\Gamma(\alpha)} \int_{0}^{\vartheta} \frac{\mathbf{c}(\zeta) \mathbf{s}(\zeta)}{\zeta^{1 - \alpha} (\vartheta - \zeta)^{1 - \alpha}} \, \rd \zeta \\
            & \hphantom{\mathbf{y}(\vartheta)}
            = \overline{\mathbf{y}}(\vartheta)
            + \frac{1}{\Gamma(\alpha)} \int_{0}^{\vartheta} \frac{\mathbf{c}(\zeta) \mathbf{y}(\zeta) + \mathbf{d}(\zeta)}{(\vartheta - \zeta)^{1 - \alpha}} \, \rd \zeta
        \end{align*}
        for all $\vartheta \in (0, 1]$.
        Thus, $\mathbf{y}(\cdot)$ is a solution of original integral equation~\eqref{integral_equation_abstract}.

        Hence, in order to complete the proof of the proposition, it suffices to show that there exists a unique function $\mathbf{s}(\cdot) \in \C[0, 1]$ satisfying integral equation~\eqref{proof_proposition_linear_existence_uniqueness_arbitrary_auxiliaty_integral_equation}.
        Since $\overline{\mathbf{s}}(\cdot) \in \C[0, 1]$, this fact can be proved by essentially repeating the arguments from~\cite[Proposition~4.1]{Gomoyunov_2020_FCAA}.
    \proofend

    Let us denote by $E_{\alpha, \alpha}$ the two-parametric Mittag--Leffler function, whose properties can be found in, e.g.,~\cite[Chapter~4]{Gorenflo_Kilbas_Mainardi_Rogosin_2014}.
    \begin{prop} \label{proposition_Gronwall}
        Let $\overline{\mathbf{y}}(\cdot) \in \C^{1 - \alpha} (0, 1]$ and $c \geq 0$ be given.
        Suppose that a function $\mathbf{y}(\cdot) \in \C^{1 - \alpha} (0, 1]$ satisfies the inequality
        \begin{equation*}
            \mathbf{y}(\vartheta)
            \leq \overline{\mathbf{y}}(\vartheta)
            + \frac{c}{\Gamma(\alpha)} \int_{0}^{\vartheta} \frac{\mathbf{y}(\zeta)}{(\vartheta - \zeta)^{1 - \alpha}} \, \rd \zeta,
            \quad \vartheta \in (0, 1].
        \end{equation*}
        Then, it holds that
        \begin{equation*}
            \mathbf{y}(\vartheta)
            \leq \overline{\mathbf{y}}(\vartheta)
            + c \int_{0}^{\vartheta} \frac{E_{\alpha, \alpha} (c (\vartheta - \zeta)^\alpha) \overline{\mathbf{y}}(\zeta)}{(\vartheta - \zeta)^{1 - \alpha}} \, \rd \zeta,
            \quad \vartheta \in (0, 1].
        \end{equation*}
    \end{prop}

    This proposition can be proved by the scheme from~\cite[Lemma~7.1.1]{Henry_1981} (see also, e.g., \cite[Lemma~2.4]{Lin_Yong_2020}).
    From Proposition~\ref{proposition_Gronwall}, we derive the following result, which, in particular, correlates with~\cite[Proposition~B.1]{Bergounioux_Bourdin_2020}.
    \begin{cor} \label{corollary_Gronwall}
        Let $\overline{y} \geq 0$, a nonnegative function $\overline{\mathbf{y}}(\cdot) \in \C^{1 - \alpha}(0, 1]$, and $c \geq 0$ be given.
        If a function $\mathbf{y}(\cdot) \in \C^{1 - \alpha} (0, 1]$ satisfies the inequality
        \begin{equation*}
            \mathbf{y}(\vartheta)
            \leq \frac{\overline{y}}{\Gamma(\alpha) \vartheta^{1 - \alpha}} + \overline{\mathbf{y}}(\vartheta)
            + \frac{c}{\Gamma(\alpha)} \int_{0}^{\vartheta} \frac{\mathbf{y}(\zeta)}{(\vartheta - \zeta)^{1 - \alpha}} \, \rd \zeta,
            \quad \vartheta \in (0, 1],
        \end{equation*}
        then the estimate below is valid:
        \begin{equation} \label{corollary_Gronwall_constants_conclusion}
            \mathbf{y}(\vartheta)
            \leq \frac{E_{\alpha, \alpha}(c) \overline{y}}{\vartheta^{1 - \alpha}}
            + \overline{\mathbf{y}}(\vartheta)
            + c E_{\alpha, \alpha}(c) \int_{0}^{\vartheta} \frac{\overline{\mathbf{y}}(\zeta)}{(\vartheta - \zeta)^{1 - \alpha}} \, \rd \zeta,
            \quad \vartheta \in (0, 1].
        \end{equation}
    \end{cor}
    \proof
        Applying Proposition~\ref{proposition_Gronwall}, we get
        \begin{align*}
            & \mathbf{y}(\vartheta)
            \leq \frac{\overline{y}}{\Gamma(\alpha)} \biggl( \frac{1}{\vartheta^{1 - \alpha}}
            + c \int_{0}^{\vartheta} \frac{E_{\alpha, \alpha} (c (\vartheta - \zeta)^\alpha)}{\zeta^{1 - \alpha} (\vartheta - \zeta)^{1 - \alpha}} \, \rd \zeta \biggr) \\
            & \hphantom{\mathbf{y}(\vartheta)}
            + \overline{\mathbf{y}}(\vartheta)
            + c \int_{0}^{\vartheta} \frac{E_{\alpha, \alpha} (c (\vartheta - \zeta)^\alpha) \overline{\mathbf{y}}(\zeta)}{(\vartheta - \zeta)^{1 - \alpha}} \, \rd \zeta
        \end{align*}
        for every $\vartheta \in (0, 1]$.
        Hence, using the equality (see, e.g., formulas (4.4.5) for $\beta = \mu = \alpha$ and (4.2.3) for $\beta = \alpha$ in \cite{Gorenflo_Kilbas_Mainardi_Rogosin_2014} and also \cite[Theorem~6.2]{Becker_2016})
        \begin{equation*}
            E_{\alpha, \alpha} (c \vartheta^\alpha)
            = \frac{1}{\Gamma(\alpha)}
            + \frac{c \vartheta^{1 - \alpha}}{\Gamma(\alpha)} \int_{0}^{\vartheta} \frac{E_{\alpha, \alpha} (c (\vartheta - \zeta)^\alpha)}
            {\zeta^{1 - \alpha} (\vartheta - \zeta)^{1 - \alpha}} \, \rd \zeta
        \end{equation*}
        and noting that $E_{\alpha, \alpha}(c \vartheta^{1 - \alpha}) \leq E_{\alpha, \alpha}(c)$, we obtain~\eqref{corollary_Gronwall_constants_conclusion}.
    \proofend

\section{Properties of solution of auxiliary integral equation}
\label{section_properties_of_solution}

\setcounter{section}{7}
\setcounter{equation}{0}\setcounter{theorem}{0}

    Fix $(t, w(\cdot)) \in G$ and $\ell \in \mathbb{R}$ and put $\lambda^{(\ell)}(\cdot) \triangleq \lambda^{(\ell)}(\cdot \mid t, w(\cdot))$ (see~\eqref{x^f}).
    For every $\tau \in [t, T)$, let $\mathbf{x}(\cdot \mid \tau) \triangleq \mathbf{x}(\cdot \mid \tau, \lambda^{(\ell)}_\tau(\cdot))$ be the solution of auxiliary integral equation~\eqref{y_integral_equation} corresponding to the pair $(\tau, \lambda^{(\ell)}_\tau(\cdot))$.

    \begin{lem} \label{lemma_mathbf_x_properties}
        The following statements are valid:

        \smallskip

        (i)
            There exists $R \geq 0$ such that, for every $\tau \in [t, T)$ and every $\vartheta \in [0, 1]$, the inequality $|\mathbf{x}(\vartheta \mid \tau)| \leq R$ takes place.

        \smallskip

        (ii)
            It holds that $\mathbf{x}(\vartheta \mid \tau) \to \mathbf{x}(\vartheta \mid t)$ as $\tau \to t^+$ uniformly in $\vartheta \in [0, 1]$.

        \smallskip

        (iii)
            For any $\eta \in (0, T - t)$, there exists $L \geq 0$ such that, for every $\tau \in (t, T - \eta]$ and every $\vartheta \in (0, 1]$, the inequality below is fulfilled:
            \begin{equation} \label{item_iii_inequality}
                \biggl| \frac{\mathbf{x}(\vartheta \mid \tau) - \mathbf{x}(\vartheta \mid t)}{\tau - t} \biggr|
                \leq \frac{L}{\vartheta^{1 - \alpha}}.
            \end{equation}
    \end{lem}
    \proof
        Consider the mapping $[t, T] \ni \tau \mapsto x(\cdot \mid \tau, \lambda^{(\ell)}_\tau(\cdot)) \in \C[0, T]$ (see Section~\ref{section_Statement_of_the_problem}), where we put formally $x(\cdot \mid T, \lambda^{(\ell)}_T(\cdot)) \triangleq \lambda^{(\ell)}(\cdot)$.
        Provided that the space $\C[0, T]$ is endowed with the uniform norm, this mapping is continuous (see the proof of~\cite[Theorem~1]{Gomoyunov_Lukoyanov_2021} and also~\cite[Theorem~2]{Gomoyunov_2020_DE}).
        In particular, the functions $x(\cdot \mid \tau, \lambda^{(\ell)}_\tau(\cdot))$ are bounded uniformly in $\tau \in [t, T]$ and $x(\xi \mid \tau, \lambda^{(\ell)}_\tau(\cdot)) \to x(\xi \mid t, w(\cdot))$ as $\tau \to t^+$ uniformly in $\xi \in [0, T]$.
        Hence, due to Lemma~\ref{lemma_mathbf_x_x}, we obtain items (i) and (ii).

        Further, recalling that the function $f$ is continuously differentiable, by the number $R$ from item (i), we can find $L_\ast \geq 0$ such that
        \begin{equation*}
            |f(\tau, \mathbf{x}) - f(\tau^\prime, \mathbf{x}^\prime)|
            \leq L_\ast (|\tau - \tau^\prime| + |\mathbf{x} - \mathbf{x}^\prime|)
        \end{equation*}
        for any $\tau$, $\tau^\prime \in [0, T]$ and any $\mathbf{x}$, $\mathbf{x}^\prime \in [- R, R]$.
        Then, given $\eta \in (0, T - t)$, we choose $\overline{L}$ according to Lemma~\ref{lemma_overline_mathbf_Lipschitz}, take $\gamma$ from condition~\eqref{sublinear_growth}, and put
        \begin{equation*}
            L^\ast
            \triangleq \Gamma(\alpha) \overline{L} + \frac{T^\alpha L_\ast}{\alpha} + \frac{\gamma (1 + R)}{\eta^{1 - \alpha}},
            \quad L
            \triangleq E_{\alpha, \alpha} (T^\alpha L_\ast) L^\ast.
        \end{equation*}

        Let $\tau \in (t, T - \eta]$ be fixed.
        Due to~\eqref{y_integral_equation}, for every $\vartheta \in (0, 1]$, we have
        \begin{align*}
            & |\mathbf{x}(\vartheta \mid \tau) - \mathbf{x}(\vartheta \mid t)|
            \leq |\overline{\mathbf{x}}(\vartheta \mid \tau) - \overline{\mathbf{x}}(\vartheta \mid t)| \\
            & \ \ + \frac{(T - \tau)^\alpha}{\Gamma(\alpha)} \int_{0}^{\vartheta}
            \frac{|f(\tau + \zeta (T - \tau), \mathbf{x}(\zeta \mid \tau)) - f(t + \zeta (T - t), \mathbf{x}(\zeta \mid t))|}
            {(\vartheta - \zeta)^{1 - \alpha}} \, \rd \zeta \\
            & \ \ + \frac{(T - t)^\alpha - (T - \tau)^\alpha}{\Gamma(\alpha)} \int_{0}^{\vartheta}
            \frac{|f(t + \zeta (T - t), \mathbf{x}(\zeta \mid t))|}{(\vartheta - \zeta)^{1 - \alpha}} \, \rd \zeta \\
            & \ \ \leq \frac{L^\ast (\tau - t)}{\Gamma(\alpha) \vartheta^{1 - \alpha}}
            + \frac{T^\alpha L_\ast}{\Gamma(\alpha)} \int_{0}^{\vartheta}
            \frac{|\mathbf{x}(\zeta \mid \tau) - \mathbf{x}(\zeta \mid t)|}{(\vartheta - \zeta)^{1 - \alpha}} \, \rd \zeta.
        \end{align*}
        Hence, applying Corollary~\ref{corollary_Gronwall}, we conclude that inequality~\eqref{item_iii_inequality} is valid for all $\vartheta \in (0, 1]$.
        The lemma is proved.
    \proofend

    By the functions $a(\cdot)$ and $b(\cdot)$ from~\eqref{a_b}, for every $\vartheta \in [0, 1]$, define
    \begin{equation} \label{mathbf_a_b}
        \mathbf{a}(\vartheta)
        \triangleq (T - t)^\alpha a(t + \vartheta (T - t)),
        \quad \mathbf{b}(\vartheta)
        \triangleq (T - t)^\alpha b(t + \vartheta (T - t)).
    \end{equation}
    Note that, since $a(\cdot)$, $b(\cdot) \in \C[t, T]$ due to the assumptions imposed on the function $f$, we have $\mathbf{a}(\cdot)$, $\mathbf{b}(\cdot) \in \C[0, 1]$.
    Consider the integral equations
    \begin{equation} \label{z^1_integral_equation}
        \mathbf{p}(\vartheta)
        = \overline{\mathbf{p}}(\vartheta) + \frac{1}{\Gamma(\alpha)}
        \int_{0}^{\vartheta} \frac{\mathbf{a}(\zeta) \mathbf{p}(\zeta) + \mathbf{b}(\zeta)}{(\vartheta - \zeta)^{1 - \alpha}} \, \rd \zeta,
        \quad \vartheta \in (0, 1],
    \end{equation}
    and
    \begin{equation} \label{z^2_integral_equation}
        \mathbf{q}(\vartheta)
        = \overline{\mathbf{q}}(\vartheta)
        + \frac{1}{\Gamma(\alpha)} \int_{0}^{\vartheta} \frac{\mathbf{a}(\zeta) \mathbf{q}(\zeta)}{(\vartheta - \zeta)^{1 - \alpha}} \, \rd \zeta,
        \quad \vartheta \in (0, 1],
    \end{equation}
    where $\overline{\mathbf{p}}(\cdot)$ and $\overline{\mathbf{q}}(\cdot)$ are the functions from~\eqref{overline_mathbf_p_q}.
    Taking into account that $\overline{\mathbf{p}}(\cdot)$, $\overline{\mathbf{q}}(\cdot) \in \C^{1 - \alpha} (0, 1]$ by Lemma~\ref{lemma_inclusions} and then applying~Proposition~\ref{proposition_linear_existence_uniqueness_arbitrary}, we conclude that integral equations~\eqref{z^1_integral_equation} and~\eqref{z^2_integral_equation} have unique solutions $\mathbf{p}(\cdot) \triangleq \mathbf{p}(\cdot \mid t, w(\cdot))$ and $\mathbf{q}(\cdot) \triangleq \mathbf{q}(\cdot \mid t, w(\cdot))$ from $\C^{1 - \alpha} (0, 1]$, respectively.

    \begin{lem} \label{lemma_mathbf_x_derivatives}
        The following limit relation holds:
        \begin{equation} \label{lemma_mathbf_x_derivatives_main}
            \lim_{\tau \to t^+} \frac{\mathbf{x}(1 \mid \tau) - \mathbf{x}(1 \mid t)}{\tau - t}
            = \mathbf{p}(1) + \mathbf{q}(1) \ell.
        \end{equation}
    \end{lem}
    \proof
        Fix $\eta \in (0, T - t)$ and choose numbers $R$ and $L$ according to items (i) and (iii) of Lemma~\ref{lemma_mathbf_x_properties}.
        Denote $\Omega \triangleq [0, 1] \times [- R, R] \times [t, T - \eta]$ and consider the continuously differentiable function
        \begin{equation*}
            \omega(\vartheta, \mathbf{x}, \tau)
            \triangleq (T - \tau)^\alpha f(\tau + \vartheta (T - \tau), \mathbf{x}),
            \quad (\vartheta, \mathbf{x}, \tau) \in \Omega.
        \end{equation*}
        Let $R_\ast \geq 0$ be such that
        \begin{equation*}
            \biggl| \frac{\partial \omega}{\partial \mathbf{x}} (\vartheta, \mathbf{x}, \tau) \biggr|
            \leq R_\ast,
            \quad (\vartheta, \mathbf{x}, \tau) \in \Omega.
        \end{equation*}

        Let $\varepsilon > 0$ be given.
        Applying Lemma~\ref{lemma_overline_mathbf_x_derivatives} and recalling notation~\eqref{overline_mathbf_z}, we can find $\delta_1 \in (0, T - \eta - t]$ such that
        \begin{equation*}
            \overline{\mathbf{z}}(1 \mid \tau)
            + R_\ast E_{\alpha, \alpha} (R_\ast) \int_{0}^{1} \frac{\overline{\mathbf{z}}(\zeta \mid \tau)}{(1 - \zeta)^{1 - \alpha}} \, \rd \zeta
            \leq \frac{\varepsilon}{2},
            \quad \tau \in (t, t + \delta_1].
        \end{equation*}
        Let us choose $\kappa > 0$ from the condition
        \begin{equation*}
            \kappa E_{\alpha, \alpha} (R_\ast) (L + 1) B(\alpha, \alpha)
            \leq \frac{\varepsilon}{2}
        \end{equation*}
        and take $\nu > 0$ such that, for any $(\vartheta, \mathbf{x}, \tau)$, $(\vartheta, \mathbf{x}^\prime, \tau^\prime) \in \Omega$, if $|\mathbf{x} - \mathbf{x}^\prime| \leq \nu$ and $|\tau - \tau^\prime| \leq \nu$, then
        \begin{equation*}
            \biggl| \frac{\partial \omega}{\partial \mathbf{x}} (\vartheta, \mathbf{x}, \tau)
            - \frac{\partial \omega}{\partial \mathbf{x}} (\vartheta, \mathbf{x}^\prime, \tau^\prime) \biggr|
            \leq \kappa,
            \quad \biggl| \frac{\partial \omega}{\partial \tau} (\vartheta, \mathbf{x}, \tau)
            - \frac{\partial \omega}{\partial \tau} (\vartheta, \mathbf{x}^\prime, \tau^\prime) \biggr|
            \leq \kappa.
        \end{equation*}
        Further, due to item (ii) of Lemma~\ref{lemma_mathbf_x_properties}, there exists $\delta_2 \in (0, \delta_1]$ such that
        \begin{equation*}
            |\mathbf{x}(\vartheta \mid \tau) - \mathbf{x}(\vartheta \mid t)|
            \leq \nu,
            \quad \vartheta \in [0, 1], \quad \tau \in [t, t + \delta_2].
        \end{equation*}
        Let us define $\delta \triangleq \min\{\nu, \delta_2\}$.

        For every $\tau \in (t, T - \eta]$ and every $\vartheta \in (0, 1]$, we denote
        \begin{equation*}
            \mathbf{r}(\vartheta)
            \triangleq \mathbf{p}(\vartheta) + \mathbf{q}(\vartheta) \ell,
            \quad \mathbf{z}(\vartheta \mid \tau)
            \triangleq \biggl|\frac{\mathbf{x}(\vartheta \mid \tau) - \mathbf{x}(\vartheta \mid t)}{\tau - t} - \mathbf{r}(\vartheta) \biggr|
        \end{equation*}
        for brevity.
        Thus, in order to prove relation~\eqref{lemma_mathbf_x_derivatives_main}, it suffices to verify that the inequality $\mathbf{z}(1 \mid \tau) \leq \varepsilon$ is valid for all $\tau \in (t, t + \delta]$.

        Let $\tau \in (t, t + \delta]$ be given.
        For any $\vartheta \in (0, 1]$, based on~\eqref{z^1_integral_equation} and~\eqref{z^2_integral_equation} and taking~\eqref{a_b}, \eqref{lemma_mathbf_x_x_main}, and~\eqref{mathbf_a_b} into account, we derive
        \begin{align*}
            & \mathbf{r}(\vartheta)
            = \overline{\mathbf{p}}(\vartheta) + \overline{\mathbf{q}}(\vartheta) \ell \\
            & \hphantom{\mathbf{r}(\vartheta)}
            + \frac{1}{\Gamma(\alpha)} \int_{0}^{\vartheta} \frac{\frac{\partial \omega}{\partial \mathbf{x}} (\zeta, \mathbf{x}(\zeta \mid t), t) \mathbf{r}(\zeta)
            + \frac{\partial \omega}{\partial \tau} (\zeta, \mathbf{x}(\zeta \mid t), t)}{(\vartheta - \zeta)^{1 - \alpha}} \, \rd \zeta,
        \end{align*}
        and, consequently, in accordance with~\eqref{y_integral_equation}, we obtain
        \begin{align}
            & \mathbf{z}(\vartheta \mid \tau)
            \leq \overline{\mathbf{z}}(\vartheta \mid \tau)
            + \frac{1}{\Gamma(\alpha)} \biggl| \int_{0}^{\vartheta} \frac{\omega(\zeta, \mathbf{x}(\zeta \mid \tau), \tau)
            - \omega(\zeta, \mathbf{x}(\zeta \mid t), t)}{(\tau - t) (\vartheta - \zeta)^{1 - \alpha}} \, \rd \zeta \nonumber \\
            & \hphantom{\mathbf{z}(\vartheta \mid \tau)} - \int_{0}^{\vartheta} \frac{\frac{\partial \omega}{\partial \mathbf{x}} (\zeta, \mathbf{x}(\zeta \mid t), t) \mathbf{r}(\zeta)
            + \frac{\partial \omega}{\partial \tau} (\zeta, \mathbf{x}(\zeta \mid t), t)}{(\vartheta - \zeta)^{1 - \alpha}} \, \rd \zeta \biggr|.
            \label{proof_mathbf_z_first}
        \end{align}

        Fix $\zeta \in (0, \vartheta]$.
        By the mean value theorem, on the segment connecting the points $(\mathbf{x}(\zeta \mid t), t)$ and $(\mathbf{x}(\zeta \mid \tau), \tau)$, there is a point $(\mathbf{x}^\prime, \tau^\prime)$ for which
        \begin{align*}
            & \omega(\zeta, \mathbf{x}(\zeta \mid \tau), \tau) - \omega(\zeta, \mathbf{x}(\zeta \mid t), t) \\
            & \quad = \frac{\partial \omega}{\partial \mathbf{x}} (\zeta, \mathbf{x}^\prime, \tau^\prime) (\mathbf{x}(\zeta \mid \tau) - \mathbf{x}(\zeta \mid t))
            + \frac{\partial \omega}{\partial \tau} (\zeta, \mathbf{x}^\prime, \tau^\prime) (\tau - t).
        \end{align*}
        Note that $(\zeta, \mathbf{x}^\prime, \tau^\prime) \in \Omega$ and, in addition, $|\mathbf{x}(\zeta \mid t) - \mathbf{x}^\prime| \leq \nu$ and $|t - \tau^\prime| \leq \nu$.
        Hence, in view of the choice of $\nu$, $R_\ast$, and $L$, we have
        \begin{align}
            & \biggl| \frac{\omega(\zeta, \mathbf{x}(\zeta \mid \tau), \tau) - \omega(\zeta, \mathbf{x}(\zeta \mid t), t)}{\tau - t} \nonumber \\
            & \quad - \frac{\partial \omega}{\partial \mathbf{x}} (\zeta, \mathbf{x}(\zeta \mid t), t) \mathbf{r}(\zeta)
            - \frac{\partial \omega}{\partial \tau} (\zeta, \mathbf{x}(\zeta \mid t), t) \biggr| \nonumber \\
            & \quad \leq \kappa \biggl( \biggl| \frac{\mathbf{x}(\zeta \mid \tau) - \mathbf{x}(\zeta \mid t)}{\tau - t} \biggr| + 1 \biggr)
            + R_\ast \mathbf{z}(\zeta \mid \tau) \nonumber \\
            & \quad \leq \frac{\kappa (L + 1)}{\zeta^{1 - \alpha}} + R_\ast \mathbf{z}(\zeta \mid \tau).
            \label{proof_mean_value_theorem}
        \end{align}

        From~\eqref{proof_mathbf_z_first} and~\eqref{proof_mean_value_theorem}, using~\eqref{formula_beta}, we derive
        \begin{equation*}
            \mathbf{z}(\vartheta \mid \tau)
            \leq \overline{\mathbf{z}}(\vartheta \mid \tau)
            + \frac{\kappa (L + 1) B(\alpha, \alpha)}{\Gamma(\alpha) \vartheta^{1 - \alpha}}
            + \frac{R_\ast}{\Gamma(\alpha)} \int_{0}^{\vartheta} \frac{\mathbf{z}(\zeta \mid \tau)}{(\vartheta - \zeta)^{1 - \alpha}} \, \rd \zeta
        \end{equation*}
        for all $\vartheta \in (0, 1]$.
        Taking into account that $\overline{\mathbf{z}}(\cdot \mid \tau)$, $\mathbf{z}(\cdot \mid \tau) \in \C^{1 - \alpha} (0, 1]$, due to Corollary~\ref{corollary_Gronwall} and the choice of $\kappa$ and $\delta_1$, we get
        \begin{equation*}
            \mathbf{z}(1 \mid \tau)
            \leq \frac{\varepsilon}{2}
            + \overline{\mathbf{z}}(1 \mid \tau)
            + R_\ast E_{\alpha, \alpha} (R_\ast) \int_{0}^{1} \frac{\overline{\mathbf{z}}(\zeta \mid \tau)}{(1 - \zeta)^{1 - \alpha}} \, \rd \zeta
            \leq \varepsilon.
        \end{equation*}
        The proof is complete.
    \proofend

\section{Proof of Theorem~\ref{theorem}}
\label{section_proof}

\setcounter{section}{8}
\setcounter{equation}{0}\setcounter{theorem}{0}

    Let $(t, w(\cdot)) \in G$ be fixed.
    First of all, note that, in view of~\eqref{overline_mathbf_p_q} and~\eqref{mathbf_a_b}, similarly to the proof of Lemma~\ref{lemma_mathbf_x_x}, it can be verified that integral equations~\eqref{integral_equation_p} and~\eqref{integral_equation_q} have unique solutions $p(\cdot) \triangleq p(\cdot \mid t, w(\cdot))$ and $q(\cdot) \triangleq q(\cdot \mid t, w(\cdot))$, respectively, which are given by
    \begin{equation*}
        p(\tau)
        = \mathbf{p} \biggl( \frac{\tau - t}{T - t} \biggr),
        \quad q(\tau)
        = \mathbf{q} \biggl( \frac{\tau - t}{T - t} \biggr),
        \quad \tau \in (t, T],
    \end{equation*}
    where $\mathbf{p}(\cdot) \triangleq \mathbf{p}(\cdot \mid t, w(\cdot))$ and $\mathbf{q}(\cdot) \triangleq \mathbf{q}(\cdot \mid t, w(\cdot))$ are the solutions of integral equations~\eqref{z^1_integral_equation} and~\eqref{z^2_integral_equation}, respectively.
    In particular, we obtain $p(T) = \mathbf{p}(1)$ and $q(T) = \mathbf{q}(1)$, and, hence, it follows from equality~\eqref{varphi_y} and Lemma~\ref{lemma_mathbf_x_derivatives} that
    \begin{equation} \label{theorem_proof_fixed_f}
        \lim_{\tau \to t^+} \frac{\rho(\tau, \lambda^{(\ell)}_\tau(\cdot)) - \rho(t, w(\cdot))}{\tau - t}
        = p(T) + q(T) \ell
    \end{equation}
    for all $\ell \in \mathbb{R}$, where $\lambda^{(\ell)}(\cdot) \triangleq \lambda^{(\ell)}(\cdot \mid t, w(\cdot))$ (see~\eqref{x^f}).

    According to~\eqref{ci-differentiability}, in order to complete the proof of the theorem, we need to take $\lambda(\cdot) \in \Lambda(t, w(\cdot))$ and show that
    \begin{equation} \label{theorem_proof_desired}
        \lim_{\tau \to t^+}
        \biggl| \frac{\rho(\tau, \lambda_\tau(\cdot)) - \rho(t, w(\cdot))}{\tau - t} - p(T) - q(T) \ell(\tau) \biggr|
        = 0,
    \end{equation}
    where
    \begin{equation} \label{theorem_proof_ell_tau}
        \ell(\tau)
        \triangleq \frac{1}{\tau - t} \int_{t}^{\tau} (^C D^\alpha_{0 +} \lambda)(\xi) \, \rd \xi,
        \quad \tau \in (t, T).
    \end{equation}

    Let us fix $\eta \in (0, T - t)$ and choose $M \geq 0$ such that $|(^C D^\alpha_{0 +} \lambda)(\xi)| \leq M$ for a.e. $\xi \in [t, T]$.
    In the proof of~\cite[Lemma~1]{Gomoyunov_Lukoyanov_2021}, it is established that the functional $\rho$ satisfies a certain local Lipschitz continuity condition, and, therefore, by~\cite[Assertion~3]{Gomoyunov_Lukoyanov_2021}, we can find $\mu_1 > 0$ and $\mu_2 > 0$ such that, for any $\lambda^\prime(\cdot)$, $\lambda^{\prime \prime}(\cdot) \in \Lambda(t, w(\cdot))$, if $|(^C D^\alpha_{0 +} \lambda^\prime)(\xi)| \leq M$ and $|(^C D^\alpha_{0 +} \lambda^{\prime \prime})(\xi)| \leq M$ for a.e. $\xi \in [t, T]$, then the estimate below holds for all $\tau \in (t, T - \eta]$:
    \begin{align}
        & | \rho(\tau, \lambda^\prime_\tau(\cdot)) - \rho(\tau, \lambda^{\prime \prime}_\tau(\cdot)) | \nonumber \\
        & \quad \leq \mu_1 \biggl| \int_{t}^{\tau}
        \bigl( (^C D^\alpha_{0 +} \lambda^\prime)(\xi) - (^C D^\alpha_{0 +} \lambda^{\prime \prime})(\xi) \bigr) \, \rd \xi \biggr|
        + \mu_2 (\tau - t)^{\alpha + 1}.
        \label{theorem_mu_1_mu_2}
    \end{align}

    In particular, for every $\tau \in (t, T - \eta]$, using definition~\eqref{theorem_proof_ell_tau} of $\ell(\tau)$ and noting that $|(^C D^\alpha_{0 +} \lambda^{(\ell(\tau))})(\xi)| = |\ell(\tau)| \leq M$ for all $\xi \in (t, T)$, we get
    \begin{equation} \label{theorem_proof_part_1}
        |\rho(\tau, \lambda_\tau(\cdot)) - \rho(\tau, \lambda^{(\ell(\tau))}_\tau(\cdot)) |
        \leq \mu_2 (\tau - t)^{\alpha + 1}.
    \end{equation}

    Now, let us prove that
    \begin{equation} \label{theorem_proof_part_2}
        \lim_{\tau \to t^+}
        \biggl| \frac{\rho(\tau, \lambda^{(\ell(\tau))}_\tau(\cdot)) - \rho(t, w(\cdot))}{\tau - t} - p(T) - q(T) \ell(\tau) \biggr|
        = 0.
    \end{equation}
    Arguing by contradiction, assume that there are a number $\varepsilon_\ast > 0$ and a sequence $\{\tau_i\}_{i \in \mathbb{N}} \subset (t, T - \eta]$ converging to $t$ as $i \to \infty$ such that
    \begin{equation} \label{theorem_proof_assumption}
        \biggl| \frac{\rho(\tau_i, \lambda^{(\ell(\tau_i))}_{\tau_i}(\cdot)) - \rho(t, w(\cdot))}{\tau_i - t} - p(T) - q(T) \ell(\tau_i) \biggr|
        \geq \varepsilon_\ast,
        \quad i \in \mathbb{N}.
    \end{equation}
    Since $|\ell(\tau_i)| \leq M$ for all $i \in \mathbb{N}$, we can suppose that there exists $\ell_\ast \in \mathbb{R}$ such that $|\ell_\ast| \leq M$  and $\ell(\tau_i) \to \ell_\ast$ as $i \to \infty$.
    Then, due to~\eqref{theorem_mu_1_mu_2}, we have
    \begin{align}
        & \biggl| \frac{\rho(\tau_i, \lambda^{(\ell(\tau_i))}_{\tau_i}(\cdot)) - \rho(t, w(\cdot))}{\tau_i - t}
        - p(T) - q(T) \ell(\tau_i) \biggr| \nonumber \\
        & \quad \leq \biggl| \frac{\rho(\tau_i, \lambda^{(\ell_\ast)}_{\tau_i}(\cdot)) - \rho(t, w(\cdot))}{\tau_i - t}
        - p(T) - q(T) \ell_\ast \biggr| \nonumber \\
        & \quad + (\mu_1 + |q(T)|) |\ell(\tau_i) - \ell_\ast| + \mu_2 (\tau_i - t)^\alpha
        \label{theorem_proof_contradiction}
    \end{align}
    for all $i \in \mathbb{N}$.
    Hence, in accordance with~\eqref{theorem_proof_fixed_f}, we obtain that the right-hand side of inequality~\eqref{theorem_proof_contradiction} tends to $0$ as $i \to \infty$, which contradicts~\eqref{theorem_proof_assumption}.

    From~\eqref{theorem_proof_part_1} and~\eqref{theorem_proof_part_2}, we derive~\eqref{theorem_proof_desired} and complete the proof.

\section{Conclusion}
\label{section_conclusion}

    In this paper, we have considered the functional $\rho$ (see~\eqref{rho_introduction} and~\eqref{rho}) that associates initial data $(t, w(\cdot)) \in G$ with the endpoint $x(T \mid t, w(\cdot))$ of the solution of Cauchy problem~\eqref{differential_equation_t}, \eqref{initial_condition_t}.
    We have proved that this functional is $ci$-differentiable of order $\alpha$ and that the corresponding derivatives can be found by solving integral equations~\eqref{integral_equation_p} and~\eqref{integral_equation_q}.

    In particular, the results of the paper illustrate that the notion of $ci$-differentiability of order $\alpha$ (see~\eqref{ci-differentiability}), which at first glance seems rather unusual, is actually quite natural and reflects the specifics of dealing with functionals defined in terms of solutions of fractional differential equations.
    Accordingly, these results, which are of interest in themselves, are also expected to find applications in the study of properties of the value functionals of optimal control problems and differential games for dynamical systems described by such equations.

\section*{Acknowledgements}

    This work was supported by RSF, project no. 19-11-00105.

\appendix

\section{Example}
\label{appendix}

    This section presents~\cite{Antonov_2020} an example of a pair $(t, w(\cdot)) \in G$ such that, for the corresponding function $\overline{\mathbf{p}}(\cdot)$ defined by~\eqref{overline_p_q} and~\eqref{overline_mathbf_p_q}, the expression $\vartheta^{1 - \alpha} \overline{\mathbf{p}}(\vartheta)$ does not have a limit as $\vartheta \to 0^+$.

    Let us suppose that $T = 2$ and take $t = 1$.
    Fix $\vartheta_\ast \in (0, 1]$ such that
    \begin{equation*}
        \int_{0}^{\frac{1}{\vartheta_\ast}} \frac{\rd u}{(1 + u)^{2 - \alpha}}
        > 2 \int_{\frac{1}{\vartheta_\ast}}^{\infty} \frac{\rd u}{(1 + u)^{2 - \alpha}}.
    \end{equation*}
    Denote
    \begin{equation*}
        \vartheta_i
        \triangleq \frac{1}{i!},
        \quad \Delta_i
        \triangleq \biggl( \frac{\vartheta_{i + 1}}{\vartheta_\ast}, \frac{\vartheta_i}{\vartheta_\ast} \biggr],
        \quad i \in \mathbb{N}.
    \end{equation*}
    Consider the function $\ell \colon [0, \frac{1}{\vartheta_\ast}] \to \mathbb{R}$ defined by
    \begin{equation*}
        \ell(\xi)
        \triangleq \begin{cases}
            0, & \mbox{if } \xi \in \bigcup\limits_{i \in \mathbb{N}} \Delta_{2 i - 1} \mbox{ or } \xi = 0, \\
            1, & \mbox{if } \xi \in \bigcup\limits_{i \in \mathbb{N}} \Delta_{2 i},
          \end{cases}
    \end{equation*}
    and put
    \begin{equation*}
        w(\tau)
        \triangleq \frac{1}{\Gamma(\alpha)} \int_{0}^{\tau} \frac{\ell(1 - \xi)}{(\tau - \xi)^{1 - \alpha}} \, \rd \xi,
        \quad \tau \in [0, 1].
    \end{equation*}
    We have $w(\cdot) \in \AC^\alpha [0, 1]$ and $(^C D^\alpha_{0 +} w)(\xi) = \ell(1 - \xi)$ for a.e. $\xi \in [0, 1]$.
    In accordance with~\eqref{overline_p_q} and~\eqref{overline_mathbf_p_q}, we obtain
    \begin{equation*}
        \vartheta^{1 - \alpha} \overline{\mathbf{p}}(\vartheta)
        = - \frac{(1 - \alpha) \vartheta^{1 - \alpha} (1 - \vartheta)}{\Gamma(\alpha)}
        \int_{0}^{1} \frac{\ell(1 - \xi)}{(1 + \vartheta - \xi)^{2 - \alpha}} \, \rd \xi,
        \quad \vartheta \in (0, 1].
    \end{equation*}
    Thus, our goal is to prove that the function
    \begin{equation*}
        \mathbf{h}(\vartheta)
        \triangleq \vartheta^{1 - \alpha} \int_{0}^{1} \frac{\ell(1 - \xi)}{(1 + \vartheta - \xi)^{2 - \alpha}} \, \rd \xi
        = \int_{0}^{\frac{1}{\vartheta}} \frac{\ell(\vartheta u)}{(1 + u)^{2 - \alpha}} \, \rd u,
        \quad \vartheta \in (0, 1],
    \end{equation*}
    does not have a limit as $\vartheta \to 0^+$.

    Due to the convergence $\vartheta_i \to 0$ as $i \to \infty$ and the choice of $\vartheta_\ast$, we can find $i_\ast \in \mathbb{N}$ and $\varepsilon_\ast > 0$ such that $\vartheta_{2 i_\ast - 1} < \vartheta_\ast$ and
    \begin{equation*}
        \int_{\frac{1}{(2 i_\ast + 1) \vartheta_\ast}}^{\frac{1}{\vartheta_\ast}} \frac{\rd u}{(1 + u)^{2 - \alpha}}
        \geq 2 \int_{\frac{1}{\vartheta_\ast}}^{\infty} \frac{\rd u}{(1 + u)^{2 - \alpha}}
        + \int_{0}^{\frac{1}{2 i_\ast \vartheta_\ast}} \frac{\rd u}{(1 + u)^{2 - \alpha}}
        + \varepsilon_\ast.
    \end{equation*}
    For every $i \in \mathbb{N}$ such that $i \geq i_\ast$, we derive
    \begin{align*}
        & |\mathbf{h}(\vartheta_{2 i}) - \mathbf{h}(\vartheta_{2 i - 1})|
        \geq \int_{0}^{\frac{1}{\vartheta_\ast}} \frac{\ell(\vartheta_{2 i} u)}{(1 + u)^{2 - \alpha}} \, \rd u
        - \int_{0}^{\frac{1}{\vartheta_\ast}} \frac{\ell(\vartheta_{2 i - 1} u)}{(1 + u)^{2 - \alpha}} \, \rd u \\
        & \hphantom{|\mathbf{h}(\vartheta_{2 i}) - \mathbf{h}(\vartheta_{2 i - 1})|}
        - \int_{\frac{1}{\vartheta_\ast}}^{\frac{1}{\vartheta_{2 i}}} \frac{\ell(\vartheta_{2 i} u)}{(1 + u)^{2 - \alpha}} \, \rd u
        - \int_{\frac{1}{\vartheta_\ast}}^{\frac{1}{\vartheta_{2 i - 1}}} \frac{\ell(\vartheta_{2 i - 1} u)}{(1 + u)^{2 - \alpha}} \, \rd u.
    \end{align*}
    Let us estimate the terms from the right-hand side of this inequality.
    Note that, for any $u \in \Bigl( \frac{1}{(2 i + 1) \vartheta_\ast}, \frac{1}{\vartheta_\ast} \Bigr]$, we have $\vartheta_{2 i} u \in \Delta_{2 i}$ and, hence, $\ell(\vartheta_{2 i} u) = 1$.
    Consequently, it holds that
    \begin{equation*}
        \int_{0}^{\frac{1}{\vartheta_\ast}} \frac{\ell(\vartheta_{2 i} u)}{(1 + u)^{2 - \alpha}} \, \rd u
        \geq \int_{\frac{1}{(2 i + 1) \vartheta_\ast}}^{\frac{1}{\vartheta_\ast}} \frac{\rd u}{(1 + u)^{2 - \alpha}}.
    \end{equation*}
    On the other hand, for any $u \in \Bigl( \frac{1}{2 i \vartheta_\ast}, \frac{1}{\vartheta_\ast} \Bigr]$, since $\vartheta_{2 i - 1} u \in \Delta_{2 i - 1}$, we obtain $\ell(\vartheta_{2 i - 1} u) = 0$.
    Therefore, we get
    \begin{equation*}
        \int_{0}^{\frac{1}{\vartheta_\ast}} \frac{\ell(\vartheta_{2 i - 1} u)}{(1 + u)^{2 - \alpha}} \, \rd u
        \leq \int_{0}^{\frac{1}{2 i \vartheta_\ast}} \frac{\rd u}{(1 + u)^{2 - \alpha}}.
    \end{equation*}
    In addition, we derive
    \begin{equation*}
        \int_{\frac{1}{\vartheta_\ast}}^{\frac{1}{\vartheta_{2 i}}} \frac{\ell(\vartheta_{2 i} u)}{(1 + u)^{2 - \alpha}} \, \rd u
        + \int_{\frac{1}{\vartheta_\ast}}^{\frac{1}{\vartheta_{2 i - 1}}} \frac{\ell(\vartheta_{2 i - 1} u)}{(1 + u)^{2 - \alpha}} \, \rd u
        \leq 2 \int_{\frac{1}{\vartheta_\ast}}^{\infty} \frac{\rd u}{(1 + u)^{2 - \alpha}}.
    \end{equation*}
    As a result, based on the estimates above and the choice of $i_\ast$ and $\varepsilon_\ast$, we conclude that $|\mathbf{h}(\vartheta_{2 i}) - \mathbf{h}(\vartheta_{2 i - 1})| \geq \varepsilon_\ast$ for all $i \in \mathbb{N}$ such that $i \geq i_\ast$, which completes the proof.

 \bigskip \smallskip

 \it

 \noindent
Krasovskii Institute of Mathematics and Mechanics \\
Ural Branch of Russian Academy of Sciences \\
S. Kovalevskaya Str., Block 16 \\
Ekaterinburg -- 620108, RUSSIA  \\[4pt]
Ural Federal University \\
Mira Str., Block 19 \\
Ekaterinburg -- 620002, RUSSIA  \\[4pt]
e-mail: m.i.gomoyunov@gmail.com
\hfill Received: November 29, 2021 \\[12pt]

\end{document}